\documentclass[11pt]{article}
\usepackage{amsfonts}
\usepackage{amsmath}

\newtheorem{theorem}{Theorem}

\newtheorem{definition}[theorem]{Definition}

\newtheorem{lemma}[theorem]{Lemma}

\newcommand{\R}{{\rm I\!R}}
\newcommand{\co}{{\rm co\,}}
\newcommand{\Limsup}{\mathop{\rm Limsup\,}}
\newcommand{\T}{^{\rm T}}

\newcommand{\shiftlist}{\setlength{\itemsep}{0pt} \setlength{\parskip}{0pt} \setlength{\leftskip}{2em}}

\title{\sc An algorithm for computing an element of the Clarke generalized Jacobian of a difference of max-type functions}
\author{Ana Horta\thanks{Instituto Polit\'ecnico de Beja},  Vera Roshchina\thanks{CIMA, Universidade de \'Evora, Portugal; Ci\^encia 2008}}

\begin{document}

\maketitle

\begin{abstract}
We show that the algorithm for computing an element of the
Clarke generalized Jacobian of a max-type function proposed by
Zheng-da Huang and Guo-chun Ma in \cite{HuangMa} can be
extended to a much wider class of functions representable as a
difference of max-type functions.
\end{abstract}

\section{Introduction}\label{s:intro}

Clarke's generalized differentiation constructions are employed
in a variety of nonsmooth optimization techniques. The relevant
generalized subdifferential and Jacobian are arguably the most
common tools applied to a wealth of essentially nonsmooth
problems (see the classical works \cite{AubinFrankowska,
Clarke90, Clarke98, RockWets98}). The new applications
employing the Clarke Jacobian are still being developed; for
example, a gradient bundle method essentially based on Clarke
subgradients has recently proved successful in solving
eigenvalue optimization problems \cite{BLO01, BLO02}, and the
nonsmooth Newton's Method \cite{QiSun93} has been applied to
important classes of nonsmooth problems, such as various
nonlinear complementarity problems, stochastic optimization,
semi-infinite programming, etc. For a brief but thorough
overview of recent applications we refer the reader to
\cite{HuangMa}.

Let a function $F:\R^n\to \R^m$ be locally Lipschitz around a
point $x\in \R^n$. By $D_F\subset \R^n$ denote the set of
points on which $F$ is differentiable. For an $\bar x\in \R^n$
let
\begin{equation}\label{eq:defClarke}
\partial F(\bar x) = \co \Limsup_{x \to \bar x\atop x \in D_F}\{\nabla F(x)\},
\end{equation}
where by $\Limsup$ we denote the {\em outer set limit} (see
\cite{RockWets98}), i.e. the union of all limits of all
converging subsequences, $\co $ is the convex hull, and $\nabla
F(x)$ is the classical Jacobian of $F$ at $x\in D_F$. The set
$\partial F(\bar x) \subset \R^{n\times m}$ is the Clarke
generalized Jacobian of $F$ at $\bar x$. For a locally
Lipschitz function this set is always nonempty and bounded (see
\cite{Clarke90}).

Let a function $F:\R^n\to \R^m$ be such that each of its
components is a pointwise maximum of a finite number of smooth
functions, i.e., $F=(f_1, f_2, \dots, f_m)$,
\begin{equation}\label{eq:MaxType}
f_i(x) = \max_{j\in J_i} f_{ij}(x), \quad \forall i\in I=\{1,\dots, m\},
\end{equation}
where $J_i$, $i\in I$ are finite index sets, and $f_{ij}:\R^n
\to \R$ are continuously differentiable for all $x\in \R^n$.


It is impossible to compute the generalized Jacobian of
max-type (or a difference of max-type) function from only the
first-order information at hand (i.e. from the gradients of the
component functions). The best one can do is to use the
estimates like the bounds in \cite{DemRubConstructive} for
quasidifferentiable functions. The algorithm suggested in
\cite{HuangMa} is essentially an elegant simplification of the
method from \cite{Gao_Clarke}; both are motivated by the
observation that for some important applications, such as
Newton's method, the computation of the whole set is not
required, we only need one arbitrary element that surely
belongs to the generalized Jacobian.

Because of the max-type structure of the function, for every
given point and any direction there exists an adjacent open set
on which the max-type function is smooth, and hence the limit
of the relevant Jacobians belongs to $\partial F$. The job of
the algorithm is to carefully select the relevant gradients to
build an element from the generalized Jacobian. In addition,
the direction in \cite{HuangMa} is chosen in a way to minimize
the computation cost. The sole goal of this paper is to
demonstrate that the original algorithm can be applied to a
wider class of functions; we do not discuss the issues of
finite precision and complexity here: this has already been
addressed in \cite{HuangMa} in detail.

We discuss the algorithm in Section~\ref{s:alg} and prove its
correctness in Section~\ref{s:proof}.

\section{The Algorithm}\label{s:alg}

Let $F:\R^n\to \R^m$ be such that $F = G-H$, where both $G$ and
$H$ are max-type functions, i.e.
\begin{equation}\label{eq:01}
G = (g_1, \dots, g_m), \quad H= (h_1, \dots, h_m)
\end{equation}
with
\begin{equation}\label{eq:02}
g_i(x) = \max_{j\in J_i}g_{ij}(x),\qquad h_i(x) = \max_{k\in K_i}h_{ik}(x), \qquad i\in I=\{1,\dots,m \},
\end{equation}
where $g_{ij}, h_{ik}:\R^n\to \R$ are $C^1$ functions, and
$J_i$ and $K_i$ are finite index sets for all $i\in I$. We will
also use the notation $F=(f_1,\dots, f_m)$ with $f_i = g_i-h_i
$, $i\in I$.

For an $x\in \R^{n}$ and each $i\in I$ define the active index
sets
$$
J_{i}\left(x\right) =\left\{ j_{0}\,\bigl|\,g_{ij_{0}}(x) =\underset{%
j\in J_i}{\max }g_{ij}(x) \right\} \mbox{ and }
K_{i}(x)
=\left\{ k_{0}\,\bigl|\,h_{ik_{0}}(x) =\underset{k\in K_i}{%
\max }h_{ik}(x) \right\}.
$$
By $\nabla f$ we denote the gradient of $f:\R^n \to \R$, and by
$e_l$ we denote the $l$-th coordinate vector: $(e_l)_l = 1$,
$(e_l)_j=0$, $l\neq j$.

The Algorithm A1 is an extension of Algorithm~2.1 in
\cite{HuangMa}. The basic idea is to consider the individual
subdifferentials of each of the functions $g_i,h_i$, $i\in I$
and from each one to choose one vertex (gradient) in a way that
all the selected vertices correspond to the same direction in
which all the aforementioned functions are differentiable.
Subroutine S1 does the selection per se, while Algorithm A1 on
Step~1 iterates through the functions $g_i,h_i$, $i\in I$ and
calls to Subroutine S1 on each iteration. On Step 2 of
Algorithm A1 the gradients selected on the previous step are
used to build an element $\xi \in \partial F(x)$.

\bigskip

\noindent{\bf Algorithm A1}
\begin{itemize}
\shiftlist
\item[\bf Input: ] A point $x\in \R^{n}$, finite index sets
    $J_i, K_i$ and functions $g_{ij}$, $j\in J_i$ and
    $h_{ik}$, $k\in K_i$, $i\in I$.

\item[\bf Step 1: ] For $i\in I$ compute $T_i(x)=${\bf
    S1}$(x, J_i, \{g_{ij}\}_{j\in J_i})$, $S_i(x)=${\bf
    S1}$(x, K_i, \{h_{ik}\}_{k\in K_i})$.

\item[\bf Step 2: ] Compute
$$\xi =\left( \nabla g_{1j_{1}}(x) +\nabla h_{1k_{1}}(x) ,\cdots, \nabla g_{mj_{m}}(x) +\nabla h_{mk_{m}}(x)
    \right) \T ,
$$
where $j_i\in T_i(x)$ and  $k_i \in S_i(x)$ are chosen
arbitrarily for each $i\in I$.

\end{itemize}

\medskip

\noindent{\bf Subroutine S1}
\begin{itemize}
\shiftlist
\item[\bf Input: ] A point $x\in \R^{n}$, a finite index
    set $J$ and functions $g_{j}$, $j\in J$.
\item[\bf Step 1$\bf '$: ] Compute the active index set
    $$J(x)=\{j_0\in J\,|\, g_{j_0}(x) = \displaystyle
    \max_{j\in J}g_j(x)\},$$ let $T^0(x)=J(x)$.
\item[\bf Step 2$\bf '$: ] For $l=1,\dots, n$ let
$$
T^{l}(x)=\left\{ t_0\in T^{l-1}(x) \,\bigl|\,\nabla g_{it_0}(x) \T e_{l}=\underset{t\in T^{l-1}(x) }{\max}\nabla g_{it }(x) \T e_{l}\right\}.
$$
Output $T^{n}(x)$.
\end{itemize}

In the next section we prove the following result.
\begin{theorem}\label{thm:main}
If $F=G-H$, where $G$ and $H$ are defined by equations
\eqref{eq:01}-\eqref{eq:02}, Algorithm A1 is well defined, and
$\xi $ generated by the algorithm is an element of $ \partial
F(x) $.
\end{theorem}

\section{Proof of the correctness of Algorithm A1}\label{s:proof}

Our proof of Theorem~\ref{thm:main} is essentially along the
lines of the proof of Theorem~2.1 in \cite{HuangMa}, albeit is
a bit shorter. We need to introduce a few definitions and
technical results first.

\begin{definition}\label{def:1.1}
A continuous mapping $ f: \R ^{n}\rightarrow \R^m $ is said to
be $PC^{1}$ on an open set $ U\subset \R ^{n}$, if there exists
a finite set of $C^{1}$ functions $ f_{j}:U\rightarrow \R
^{m}$, $ j\in J$ (with $|J|<\infty$), such that for every $
x\in U$, $f(x) =f_{j}(x)$ for at least one index $ j\in J$.
\end{definition}

Recall that a directional derivative of a function $f:\R^n\to
\R$ at a point $x\in \R^n$ in the direction $y$ is the quantity
\begin{equation}\label{eq:DefDD}
f'(x;y) = \lim_{t\downarrow 0}\frac{f(x+ty)-f(x)}{t}.
\end{equation}
All $PC^1$ functions are directionally differentiable, i.e. the
limit \eqref{eq:DefDD} exists for all directions $y\in \R^n$.
The next result follows from the definition of the directional
derivative. For a detailed discussion see
 \cite[Chapter I, Corollary~3.2]{DemRubConstructive}.
\begin{lemma}\label{lem:max_dd} Let $f:\R^n\to \R$ be a pointwise
maximum of a finite number of smooth functions, i.e. for all
$x\in \R^n$
$$
f(x) ={\max_{j\in J} }f_{j}(x),
$$
where $f_{j}:
\R ^{n}\rightarrow %
\R$, $ j\in J$ are continuously differentiable, and $J$ is a
finite index set. Then for every $x\in \R^n $ the function $f$
is directionally differentiable along an arbitrary direction
$y\in \R^{n}$, and
\begin{equation}\label{eq:3.1}
f^\prime(x;y) =\max_{j\in J(x)} f^\prime _{j}( x;y) = \max_{j\in J(x)}\nabla
f_{j}(x) \T y,
\end{equation}
where $J(x)$ is the active index set:
$$
J(x) = \{j_0\in J\,|\, f_{j_0}(x) = \max_{j\in J}f_j(x)\}.
$$
\end{lemma}

The following result is proved in \cite[Lemma~2]{PangRalph96}
(also see \cite{HuangMa}).
\begin{lemma}\label{lem:JacDD} Let $F:\R^n\to \R^m$ be a $PC^1$ function in a
neighborhood of $x\in \R^n$, then
$$
\partial [F'(x;\cdot)](0)\subset \partial F(x),
$$
where $[F'(x;\cdot)]:\R^n\to \R^m$ is the vector function of
the directional derivatives of the components of $F$ at the
point $x$.
\end{lemma}

We are now in the position to prove our main result. The proof
essentially follows the ideas of the proof of Theorem~2.1 in
\cite{HuangMa}.

\bigskip

\noindent{\sc Proof of Theorem~\ref{thm:main} } For each $ i\in
I$ let $T_{i}^{0}(x):=J_{i}(x)$, $S_{i}^{0}(x):=K_{i}(x)$, and
for each $l\in \{1, \dots, n\}$ define the index subsets
$T_i^l(x)$ and $S_i^l(x)$ recursively
\begin{equation}\label{eq:3.5}
 T_{i}^{l}(x) =\left\{ t\in T_{i}^{l-1}(x)\, \left|\,\nabla g_{it}(x) \T e_{l}=\underset{j\in
T_{i}^{l-1}(x) }{\min }\nabla g_{ij}(x) \T e_{l}\right.\right\};
\end{equation}
\begin{equation}\label{eq:3.7}
 S_{i}^{l}(x) =\left\{ s\in S_{i}^{l-1}(x)\, \left|\,\nabla h_{is}(x) \T e_{l}=\underset{k\in
S_{i}^{\,l-1}(x) }{\min}\nabla h_{ik}(x) \T e_{l}\right.\right\}.
\end{equation}
It is not difficult to observe that the sets $T_i(x) =
T_i^n(x)$ and $S_i(x) = S_i^n(x)$ are precisely the sets
obtained after the execution of Step~1 of Algorithm~A1. Since
the index sets  $J_{i}(x)$ and $K_i(x)$ are nonempty and
finite, the minimal values of the scalar products in
\eqref{eq:3.5} and \eqref{eq:3.7} are attained, and hence on
every iteration of Step~2$'$ of Subroutine S1 we generate
nonempty sets. This means that after the execution of Step 1 of
A1 we end up with nonempty finite sets $T_i^n(x)$ and
$S_i^n(x)$, $i\in I$, so we can choose the corresponding
indices $j_1,\dots, j_m,k_1,\dots, k_m$ on Step~2.

Let
\begin{align*}
\Gamma := &\left\{ \nabla g_{ij}(x) -\nabla
 g_{it}(x) \,|\,j\in J_{i}(x) \setminus T_{i}(x) ,t\in T_{i}(x)
 ,i\in I\right\}\\
&\cup \left\{\nabla
 h_{ik}(x) -\nabla h_{is}(x)\,|\,k\in K_{i}(x) \setminus
 S_{i}(x) ,s\in S_{i}(x) ,i\in I\right\}.
\end{align*}
It follows from \eqref{eq:3.5}-\eqref{eq:3.7} that all the
first nonzero components of all elements in $\Gamma$ are
positive.

Let $\varepsilon $ denote a positive number smaller than the
minimum value among the first nonzero components of all
elements in $\Gamma $ and let $M$ be a positive number larger
than the maximum value among the absolute values of all
components of all elements in $\Gamma$, i.e. for each $\alpha =
\left( 0,\ldots ,0,\alpha _{k},\alpha _{k+1},\ldots ,\alpha
_{n}\right) \T \in \Gamma$ we have $| \alpha _{i}| <M,$
$i=k\ldots ,n$ and $ \alpha _{k}>\varepsilon >0$.

Let
$$\overline{y}=\left( -\lambda _{1},-\lambda _{2},\ldots
,-\lambda _{n}\right) \T ,$$ where
$\lambda _{i}>0,$ \ $\frac{\lambda _{i+1}}{\lambda _{i}}<\frac{\frac{%
\varepsilon }{M}}{1+\frac{\varepsilon }{M}},$ \ $i=1,\ldots
,n-1$. We have
\begin{align*}
\alpha \T \overline{y} & =-\underset{i=k}{\overset{n}{\sum
}}\lambda _{i}\alpha _{i}\\
 & \leq -\lambda _{k}\alpha _{k}+\underset{i=k+1}{\overset{n}{\sum }}%
\lambda _{i}| \alpha _{i}|\\
& \leq -\lambda _{k}\alpha _{k}+M\underset{i=k+1}{\overset{n}{\sum }%
}\lambda _{i}\\
&=-\lambda _{k}\alpha _{k}+M \lambda _{k}\underset{i=k+1}{%
\overset{n}{\sum }}\left( \frac{\frac{\varepsilon
}{M}}{1+\frac{\varepsilon }{M}}\right) ^{i-k}\\
& <-\lambda _{k}\alpha _{k}+M \left( \frac{\varepsilon
}{M}\right) \lambda _{k}\\
& =\lambda _{k}\left( \varepsilon -\alpha _{k}\right) \\
& <0,
\end{align*}
which means that $ \alpha \T \bar y<0$ for all $ \alpha \in
\Gamma$. Observe that the set
$$
U = \{y\,|\, \alpha\T y<0;\forall \alpha \in \Gamma\}
$$
is an open convex cone, which is nonempty since $\bar y \in U$.
Hence, we have
$$
\left( \nabla g_{ij}(x) -\nabla g_{it}(x) \right) \T y<0, \ \forall j\in J_{i}(x) \setminus T_{i}(x) ,t\in T_{i}(x)
$$
and
$$
\left( \nabla h_{ik}(x) -\nabla
h_{it}(x) \right) \T y<0, \forall k\in
K_{i}(x) \setminus S_{i}(x) ,s\in
S_{i}(x),
$$
which implies
\begin{equation}\label{eq:3.21}
\underset{j\in J_{i}(x) }{\max }\nabla
g_{it}(x) \T y=\underset{t\in T_{i}(x)
}{\max }\nabla g_{ij}(x) \T y
\end{equation}
and
\begin{equation}\label{eq:3.22}\underset{k\in K_{i}(x) }{\max }\nabla
h_{is}(x) \T y=\underset{s\in S_{i}(x)
}{\max }\nabla h_{ik}(x) \T y
\end{equation}
for all $y\in U$. It follows directly from \eqref{eq:3.5} and
\eqref{eq:3.7} that for every $i\in I$ we have
\begin{equation}\label{eq:grad_coincide} \nabla g_{i{t_1}}(x) =
\nabla g_{i{t_2}}(x)\; \forall t_1,t_2\in T_i(x);\quad \nabla h_{i{s_1}}(x) =
\nabla h_{i{s_2}}(x)\; \forall s_1,s_2\in S_i(x).
\end{equation}
Hence, for each $i\in I$ we get
\begin{equation}\label{eq:3.23a}
f^\prime _{i}\left( x;y\right) = \nabla g_{ij}(x)
 \T y-\nabla h_{ik}(x) \T y =\xi_i\T y
\end{equation}
for all $y\in U$, and any arbitrary choice of $j\in
 T_{i}(x)$ and $k\in S_{i}(x)$, i.e. $f'(x;y)$ is linear in $y$
 on $U$.

We next show that $\xi \in \partial _{B}F(x)$. Fix an arbitrary
$y\in U$. Since $U$ is an open cone, for all $t>0$ we have
$ty\in U$ together with a neighborhood; then from
\eqref{eq:3.23a} we have $\nabla [F'(x;\cdot)](ty) = \xi$, and
hence
\begin{equation}\label{eq:3.25}
\xi = \lim\limits_{t\downarrow 0} \nabla[F'\left( x;\cdot\right)](ty)
\in \partial [F'(x;\cdot)](0) \subset \partial F(x),
\end{equation}
where the last inclusion follows from  Lemma~\ref{lem:JacDD}
and the observation that $F$ is $PC^1$. The proof is complete.

\bibliographystyle{plain}	
\bibliography{references}

\end{document}